\documentclass[12pt]{amsart} 
\usepackage{amsmath,amssymb,amscd,array}

\let\ssection=\section
\renewcommand{\section}{\setcounter{equation}{0}\ssection}

\def\d{\delta}

\def\om{\omega}
\def\r{\rho}
\def\a{\alpha}
\def\b{\beta}
\def\s{\sigma}
\def\vfi{\varphi}

\def\l{\lambda}
\def\m{\mu}

\def\implies{\Rightarrow}

\newcommand{\Diff}{\mathrm{Diff}}

\newcommand{\cF}{{\mathcal{F}}}

\begin{document}

\frenchspacing

\def\d{\delta}
\def\g{\gamma}
\def\om{\omega}
\def\r{\rho}
\def\a{\alpha}
\def\b{\beta}
\def\s{\sigma}
\def\vfi{\varphi}
\def\l{\lambda}
\def\m{\mu}
\def\implies{\Rightarrow}

\oddsidemargin .1truein
\newtheorem{thm}{Theorem}[section]
\newtheorem{lem}[thm]{Lemma}
\newtheorem{cor}[thm]{Corollary}
\newtheorem{pro}[thm]{Proposition}
\newtheorem{ex}[thm]{Example}
\newtheorem{rmk}[thm]{Remark}
\newtheorem{defi}[thm]{Definition}
\title[Ternary differential operators]{The ternary invariant differential
operators acting on the spaces
of weighted densities }
\author{Sofiane BOUARROUDJ}\address{Department of Mathematical Sciences, U.A.E.
University, P.O. Box 17551, Al-Ain, U.A.E.\\
e-mail:bouarroudj.sofiane@uaeu.ac.ae}

\begin{abstract}
Over $n$-dimensional manifolds, I classify ternary differential
operators acting on the spaces of weighted densities and invariant
with respect to the Lie algebra of vector fields. For $n=1$, some of
these operators can be expressed in terms of the de Rham exterior
differential, the Poisson bracket, the Grozman operator and the
Feigin-Fuchs anti-symmetric operators; four of the operators are
new, up to dualizations and permutations.

For $n>1$, I list multidimensional conformal tranvectors, i.e.,
operators acting on the spaces of weighted densities and invariant
with respect to $\mathfrak{o}(p+1,q+1)$, where $p+q=n$. Except for
the scalar operator, these conformally invariant operators are not
invariant with respect to the whole Lie algebra of vector fields.
\end{abstract}


\thanks{I am thankful to D.~Leites and V.~Ovsienko for their interest in this work and
help.}

\maketitle

\section{Introduction}
Let $M^n$ be a manifold; denote by $\mathfrak{vect}(M)$ the Lie
algebra of vector fields on $M$ with coefficients to be specified
later (smooth, polynomial, etc.). Let $T(V)$ be the
$\mathfrak{vect}(M)$-module of tensor fields of type $V$ for a
representation $\rho: GL(n)\rightarrow GL(V)$; for a precise
description of the $\mathfrak{vect}(M)$-action, see \cite{gls, hl}.
Set $\mathfrak{vect}(n):=\mathfrak{vect}(\mathbb{R}^n)$.

All laws of physics are described by operators $D:
T(V_1)\otimes\cdots\otimes T(V_k)\rightarrow T(W)$ {\it invariant}
with respect to the action of $\mathfrak{vect}(M)$ or its
subalgebra, hence their importance (\cite{kms}). The classification
of invariant operators is a problem raised by Veblen in 1928, at the
IMC, Bologna, see \cite{ki}. Since then, many operators earlier
known in particular cases have been generalized and several new ones
were found.

In this note I will list ternary $\mathfrak{vect}(1)$-invariant
differential operators acting on the tensor fields on $\mathbb{R}$.
Clearly, on 1-dimensional manifolds, the spaces $T(V)$ are  of the
form $\cF_\l=\left\{fdx^{\lambda}\mid f\in \mathbb{R}[x]\right\}$
--- the spaces of weighted densities of weight $\l\in\mathbb{R}$.

In the unary case, Rudakov \cite{rud}, Terng \cite{te},  Kirillov
\cite{ki} have independently and by different means proved that
there are only one invariant differential operator --- the {\it
exterior differential}. In particular, for $n=1$, we have:
\[
d:\cF_0\rightarrow \cF_1, \quad \varphi\mapsto \varphi'dx.
\]
On the space of functions with compact support, or if $M$ is
compact, there is also a non-differential operator, the {\it
integral}:
\[\int: \cF_1\rightarrow \mathbb{R}\subset\cF_0, \quad
\varphi\, dx\mapsto \int \varphi\, dx\quad \text{(a non-differential
one).}\] In what follows, all functions considered are polynomial
ones and the operators are differential ones. 
How to describe non-differential operators is unclear, and
Kirillov's example of a non-differential operator different from the
integral, see \cite{ki, gls}, adds to the mystery.

Observe that (1) the following scalar operator is invariant for any
\lq\lq arity" $k$:
\begin{equation}\label{scal}
\begin{array}{ccl}
\mathcal{F}_{\lambda_1}\otimes\dots\otimes
\mathcal{F}_{\lambda_k}&\rightarrow&\mathcal{F}_{\lambda_1+\dots+\lambda_k}\\[1mm]
\varphi_1(dx)^{\lambda_1}\otimes\dots\otimes\varphi_k(dx)^{\lambda_k}&\mapsto&
\left(\prod\varphi_{i}\right)(dx)^{\lambda_1+\dots+\lambda_k}; \\
  \end{array}
\end{equation}
(2)  every invariant {\it differential} operator has constant
coefficients, see \cite{ki}.

In the binary case, the classification (for any $n$) is due to
Grozman \cite{Gz}. In particular, for $n=1$, we have:

\underline{First order operators}:
 \begin{equation*}\label{ord1}
  \begin{array}{ll}
\mathcal{F}_0\otimes\mathcal{F}_0\rightarrow\mathcal{F}_1,&
\varphi\otimes\psi\mapsto (a\varphi'\psi+b\varphi\psi')dx\text{~~for
any $a,\,b\in\mathbb{R}$},
\\[1mm]
\mathcal{F}_\lambda\otimes\mathcal{F}_\mu\rightarrow\mathcal{F}_{\lambda+\mu+1}&
\varphi(dx)^{\lambda}\otimes\psi(dx)^{\mu}\mapsto
(-\lambda\varphi\psi'+\mu\varphi'\psi)(dx)^{\lambda+\mu+1}. \\
  \end{array}
  \end{equation*}
Observe that the second of these operators is the Poisson bracket in
coordinates $x$ and $p:=dx$. In what follows I denote this operator
$\{\varphi(dx)^{\lambda}, \psi(dx)^{\mu}\}$.

\underline{Second order operators}:
  \begin{equation*}\label{ord2}
  \begin{array}{ll}
\mathcal{F}_0\otimes\mathcal{F}_\mu\rightarrow\mathcal{F}_{\mu+2}&
\varphi\otimes\psi(dx)^{\mu}\mapsto  (-\varphi'\psi'+\mu\varphi''\psi)(dx)^{\mu+2},
\\[1mm]
\mathcal{F}_\lambda\otimes\mathcal{F}_0\rightarrow\mathcal{F}_{\lambda+2}&
\varphi(dx)^{\lambda}\otimes\psi\mapsto  (-\lambda\varphi\psi''+\varphi'\psi')
(dx)^{\lambda+2}, \\[1mm]
\mathcal{F}_\lambda\otimes\mathcal{F}_{-\lambda-1}\rightarrow\mathcal{F}_{1}&
\varphi(dx)^{\lambda}\otimes\psi(dx)^{-\lambda-1}\mapsto\\
&
\big(-\lambda\varphi\psi''-(2\lambda+1)\varphi'\psi'-(\lambda+1)\varphi''\psi\big)dx. \\
  \end{array}
 \end{equation*}

\underline{Third order operators}:
\begin{equation*}\label{ord3}
 \begin{array}{ll}
 \mathcal{F}_0\otimes\mathcal{F}_0\rightarrow\mathcal{F}_3&
\varphi\otimes\psi\mapsto
(\varphi'\psi''-\varphi''\psi')(dx)^{3},\\[1mm]
\mathcal{F}_0\otimes\mathcal{F}_{-2}\rightarrow
\mathcal{F}_{1}&\varphi\otimes\psi(dx)^{-2}\mapsto
(\varphi'\psi''+3\varphi''\psi'+2\varphi'''\psi)dx, \\[1mm]
\mathcal{F}_{-2}\otimes\mathcal{F}_{0}\rightarrow\mathcal{F}_{1}&
\varphi(dx)^{-2}\otimes\psi\mapsto  (\varphi''\psi'+3\varphi'\psi''+
2\varphi\psi''')dx, \\[1mm]
\mathrm{Gz}:\mathcal{F}_{-\frac{2}{3}}\otimes\mathcal{F}_{-\frac{2}{3}}
\rightarrow\mathcal{F}_{\frac{5}{3}}&
\varphi(dx)^{-\frac{2}{3}}\otimes\psi(dx)^{-\frac{2}{3}}\mapsto\\[1mm]
&(2\varphi'''\psi+3\varphi''\psi'-3\varphi'\psi''-2\varphi\psi''')(dx)^{\frac{5}{3}}.
\end{array}
\end{equation*} Thus, for $n=1$, all invariant operators are of order $\leq 3$
and can be expressed as a composition of the de Rham differential
and the Poisson bracket, except for the {\it Grozman operator}
$\mathrm{Gz}$. This operator was the starting point for Feigin and
Fuchs \cite{ff} in their classification of all anti-symmetric
$k$-nary differential operators.

For $n>1$, I will show that there are no
$\mathfrak{vect}(n)$-invariant operators acting on the spaces of
weighted densities, except the scalar one (\ref{scal}). We can,
however, look for invariance with respect to some subalgebras of
$\mathfrak{vect}(n)$. Most interesting are the maximal simple
subalgebras.

For $n=1$,  there is the only one maximal simple subalgebra:
$\mathfrak{sl}(2)\simeq\mathfrak{o}(1,2)$. The operators invariant
with respect to it are already classified: in the unary case, these
are {\it Bol's operators} \cite{bol}; in the binary case, these are
{\it Gordon's transvectants} \cite{go}, for their application, see
\cite{bo, co, CMZ, ol, ra, za}; for the ternary case, see \cite{b2}.

For $n>1$, we can consider multidimensional analogs of transvectants
invariant with respect to various maximal simple finite dimensional
subalgebras of $\mathfrak{vect}(n)$, e.g., $\mathfrak{sl}(n+1)$ or
$\mathfrak{o}(p+1,q+1)$, where $n=p+q$ (cf., \cite{bm}). The
$\mathfrak{o}(p+1,q+1)$-invariant operators have been intensively
studied. In the unary case, all conformally invariant operators are
given by powers of the Laplacian (cf. \cite{er,gh,jv}):
\[
\Delta^{k}_{g}:\cF_\lambda\rightarrow \cF_{\lambda+\frac{2k}{n}},
\]
where $\Delta_{g}$ is the Laplacian associated with a
pseudo-Riemannian metric $g$ of signature $(p,q)$. Ovsienko and
Redoux \cite{or} classified the $\mathfrak{o}(p+1,q+1)$-invariant
binary operators. They proved that, for almost all values of the
weights, the space of conformally invariant operators is
one-dimensional.

Observe that $\mathfrak{vect}(1)\simeq\mathfrak{k}(1)$; so it is
also natural to generalize the quest for transvectants \lq\lq in the
contact direction", as in \cite{mat}, where the bilinear
differential operators acting on the spaces of weighted densities
and invariant with respect to the symplectic Lie algebra
$\mathfrak{sp}(2n+2)\subset\mathfrak{k}(2n+1)$ are classified.

\subsection{Main results} I classify all
 $\mathfrak{vect}(n)$-invariant
ternary differential operators acting on the spaces of weighted
densities. Some of these invariant operators can be expressed in
terms of unary and binary ones. Some of them generalize the
Feigin-Fuchs anti-symmetric invariant operators; four (up to
dualizations and permutations) invariant operators are new.

I also classify all conformally invariant ternary differential
operators acting on the spaces of weighted densities.
Straightforward computations and repeated arguments are omitted.

\section{Invariant operators on $\mathbb{R}$ or $S^1$}

\subsection{The Feigin-Fuchs anti-symmetric invariant operators}
In \cite{ff}, Feigin and Fuchs  classified all anti-symmetric
$k$-nary invariant differential operators on weighted densities.
Recall that their list of ternary operators consists of
$$
\begin{array}{lll}
\Delta_{\l,3}:\wedge^3 \cF_{\lambda}\rightarrow \cF_{3\l+3},& d\circ
\Delta_{-1,3}:\wedge^3 \cF_{-1}\rightarrow \cF_{1},&
\Upsilon:\wedge^3 \cF_{-\frac{5}{4}}\rightarrow \cF_{\frac{9}{4}},
\\[3mm]
\Delta_{1,3}\circ (d\otimes d\otimes d):\wedge^3 \cF_{0}\rightarrow
\cF_{6},& \Theta_{\pm}:\wedge^3 \cF_{\kappa} \rightarrow
\cF_{3\kappa+5}&\text{for $\kappa=-\frac{1}{12}(9\pm \sqrt{21})$}
\end{array}
$$
and where the operators  $\Delta_{\l,3}(\varphi, \psi, \chi)$,
$\Upsilon$, and 
$\Theta_{\pm}$ 
are defined by the
following expressions, respectively:
\[\small
\begin{array}{lcl}
&& \left |
\begin{array}{lll}
\varphi & \psi & \chi\\
\varphi' & \psi' &\chi' \\
\varphi'' & \psi'' & \chi''
\end{array}
\right |\,(dx)^{3\lambda+3},\\
&&\left (\left |
\begin{array}{lll}
\varphi & \psi & \chi\\
\varphi' & \psi' & \chi'\\
\varphi^{(5)} & \psi^{(5)} & \chi^{(5)}\end{array} \right | +
\frac{5}{2} \left |
\begin{array}{lll}
\varphi & \psi & \chi\\
\varphi'' & \psi'' & \chi''\\
\varphi^{(4)} & \psi^{(4)} & \chi^{(4)} \end{array} \right |+ 2
\left |
\begin{array}{lll}
\varphi' & \psi' & \chi'\\
\varphi'' & \psi'' & \chi''\\
\varphi^{(3)} & \psi^{(3)} & \chi^{(3)} \end{array}\right
|\right)(dx)^{\frac{9}{4}},\\[5mm]
&&\left (\left |
\begin{array}{lll}
\varphi & \psi & \chi\\
\varphi' & \psi' & \chi'\\
\varphi^{(4)} & \psi^{(4)} & \chi^{(4)}\end{array} \right |
 + 2(\mp\sqrt{21}-4) \left |
\begin{array}{lll}
\varphi & \psi & \chi\\
\varphi'' & \psi'' & \chi''\\
\varphi^{(3)} & \psi^{(3)} & \chi^{(3)} \end{array} \right |\right
)(dx)^{3\kappa+5}.
\end{array}
\]
\subsection{The list of $\mathfrak{vect}(1)$-invariant ternary
invariant differential operators}
Recall that with any operator $A:{\mathcal F}_{\lambda}\otimes
{\mathcal  F}_{\gamma}\otimes {\mathcal  F}_\tau \longrightarrow
{\mathcal  F}_{\mu}$ we can associate

(i) Three dualizations ($i=1, 2, 3$, cf. \cite{b3}):
\[\small A^{*i}:{\mathcal  F}_{\alpha}\otimes {\mathcal  F}_{\beta}\otimes {\mathcal
F}_\nu \longrightarrow {\mathcal  F}_{\delta}, \] where $(\alpha,
\beta, \nu, \delta)=(1-\mu, \gamma, \tau, 1-\lambda)$ or $(\lambda,
1-\mu, \tau, 1-\gamma)$ or $(\lambda, \gamma, 1-\mu, 1-\tau)$. These
dualizations exist due to the definition according to which
$({\mathcal F}_{\alpha})^*:={\mathcal F}_{1-\alpha}$. (Assuming that
one of these spaces consists of functions with compact support to
enable integration one can justify this definition, more correct
than a seemingly more natural $({\mathcal F}_{\alpha})^*={\mathcal
F}_{-\alpha}$.)

(ii) Permutations:

\[A^{\sigma}:{\mathcal  F}_{\alpha}\otimes {\mathcal  F}_{\beta}\otimes {\mathcal
F}_\nu \longrightarrow {\mathcal  F}_{\mu},\quad A^{\sigma}(a, b,
c)=A(\sigma^{-1}(a, b, c)),
\]
where $(\alpha, \beta, \nu)=\sigma(\lambda, \gamma, \tau)$ for any
$\sigma\in S_3$. 
\subsubsection{Order 1 and 2}
Hereafter all $\epsilon$'s are constants.

\begin{thm} (i) Every order $1$ invariant differential operator is, up to
permutations and dualizations, a multiple of one of the following
operators $$\varphi\,(dx)^{\lambda} \otimes \psi\,(dx)^{\gamma}
\otimes
\chi\,(dx)^{\tau}\mapsto(...)(dx)^{\lambda+\gamma+\tau+1},\text{~~
where $(...)=$:}
$$
$$
\begin{array}{l}
\begin{cases}
\left(\epsilon_1\, d \varphi \cdot \psi\cdot
\chi+\epsilon_2\,\varphi \cdot d \psi \cdot \chi+\epsilon_3\,\varphi
\cdot \psi\cdot d \chi\right) & \text{ for }
(\lambda, \gamma,\tau)=(0,0,0),\\[4mm]
\epsilon_1\left \{\varphi(dx)^\lambda,\psi(dx)^\gamma \right
\}\chi(dx)^{\tau}+\epsilon_2\left \{\varphi(dx)^\lambda,
\chi(dx)^\tau\right \}\psi(dx)^{\gamma} &\text{ for } \lambda\not
=0.
\end{cases}
\end{array}
$$
(ii) Every order $2$ invariant differential operator is, up to
permutations and dualizations, a multiple of one of the following
operators:
$$
\begin{array}{l}
\varphi\,(dx)^{\lambda} \otimes \psi\,(dx)^{\gamma} \otimes
\chi\,(dx)^{\tau}\mapsto\\
\begin{cases}
 \left(\epsilon_1\left \{\varphi,d \psi\right \}\chi(dx)^{\tau}+\epsilon_2\left
\{\chi (dx)^\tau,d \varphi \right \}\psi+\epsilon_3\left \{d
\psi,\chi (dx)^\tau \right \}\varphi\right)
& \text{for }(\gamma,\lambda)=(0,0),\\[3mm]
\epsilon_1\, \left \{ \left\{\varphi (dx)^{\lambda},\psi
(dx)^{\gamma}\right \},\chi (dx)^{\tau}\right
\} & \\[2mm]
+\epsilon_2\,\left \{
\left\{\varphi(dx)^{\lambda},\chi(dx)^{\tau}\right \} ,\psi
(dx)^{\gamma} \right \}\,&\text{otherwise.}
\end{cases}
\end{array}
$$
\end{thm}
\subsubsection{Order 3}
The following ternary $\mathfrak{vect}(1)$-invariant differential
operator is new, and so are its dualizations and permutations (here
$\mu=\lambda+\gamma+\tau+3$):
\[
\small \Delta_{\lambda, \gamma, \tau;3} :\cF_\lambda\otimes
\cF_\gamma\otimes \cF_\tau \rightarrow \cF_{\mu} \quad
\varphi\,(dx)^\lambda\otimes \psi\,(dx)^\gamma \otimes
\chi\,(dx)^\tau \mapsto
\sum_{i+j+l=3}\!\!\!\alpha_{i,j,l}\,\varphi^{(i)}\psi^{(j)}\chi^{(l)}\,
(dx)^{\mu},\] where the coefficients are given by
\[
\tiny{
\begin{array}{ll}
\alpha_{0,3,0}=\lambda\tau(\tau-\lambda)(1+\lambda+\tau),&
\alpha_{2,1,0}=\tau(1+\gamma+\tau) (2 \lambda +\tau +\gamma  (3
\lambda +3 \tau +4)+2),
 \\[2mm]

\alpha_{1,1,1}=(\lambda-\gamma)(\gamma-\tau)(\lambda-\tau),&
\alpha_{1,2,0}=-\tau(1+\lambda+\tau) (2 \gamma +\tau +\lambda  (3
\gamma +3 \tau +4)+2 ),
\\[2mm]

\alpha_{3,0,0}=\gamma\tau(\tau-\gamma)(1+\gamma+\tau),&
\alpha_{2,0,1}=-\gamma(1+\gamma+\tau)
(\gamma +2 \lambda +(3 \gamma +3 \lambda +4) \tau +2),\\[2mm]

\alpha_{0,0,3}=\gamma\lambda(\lambda-\gamma)(1+\gamma+\lambda),&
\alpha_{1,0,2}=\gamma(1+\gamma+\lambda)
(\gamma +2 \tau +\lambda  (3 \gamma +3 \tau +4)+2),\\[2mm]

\alpha_{0,2,1}=\lambda(1+\lambda+\tau) (2 \gamma +\lambda +(3 \gamma
+3 \lambda +4) \tau +2 ), &\alpha_{0,1,2}=-\lambda(1+\gamma+\lambda)
(\lambda +2 \tau +\gamma(3\lambda+3\tau+4)+2).
\end{array}
}
\]
The following proposition shows that $\Delta_{\lambda, \gamma,
\tau;3}$ generalizes the Feigin-Fuchs operator $\Delta_{\lambda,3}$.
\begin{pro}
(i) If $\lambda=\gamma=\tau$, then $$\Delta_{\lambda, \gamma,
\tau;3}=-\lambda (1+2 \lambda)^2 (2+3 \lambda) \Delta_{\lambda,3}.$$

(ii) If $1+\lambda+\gamma=0$, then $$\Delta_{\lambda, \gamma,
\tau;3}=(\tau-\lambda)(1+\lambda+\tau) \left \{ d \left
\{\cdot,\cdot \right \},\cdot\right \}.$$ Analogs of this equality
follow by dualizations and permutations.
\end{pro}

\begin{thm} \label{mainthm} Every order $3$ ternary invariant differential operator
is, up to permutations and dualizations, a multiple of one of the
following operators:
$$
\begin{array}{l}
\varphi\,(dx)^{\lambda} \otimes \psi\,(dx)^{\gamma} \otimes
\chi\,(dx)^{\tau}\mapsto\\
\begin{cases}
(\epsilon_1\left \{d\varphi,d\psi\right \}\chi+\epsilon_2\left
\{d\chi,d\varphi\right \}\psi+\epsilon_3\left \{d\psi,d\chi \right
\}\varphi)+ & \\[2mm]\epsilon_4 d\varphi\otimes d\psi\otimes d\chi&
\text{for }(\gamma,\lambda,\tau)=(0,0,0), \\[3mm]
\epsilon_1\, \left \{ \left \{ d \varphi,\psi(dx)^\gamma\right
\},\chi(dx)^\tau\right \}+&\\[2mm]
\epsilon_2 \, \left \{ \left \{d\varphi, \chi (dx)^\tau \right
\},\psi(dx)^\gamma\right \}
&\text{for }\lambda=0  \text{ and $\gamma^2+\tau^2\not =0$ }, \\[3mm]
\epsilon_1 \mathrm{
Gz}(\varphi\,(dx)^{-\frac{2}{3}},\psi\,(dx)^{-\frac{2}{3}})\chi(dx)^{-\frac{2}{3}}+\\[2mm]
\epsilon_2\,\mathrm{Gz}(\chi\,(dx)^{-\frac{2}{3}},
\varphi\,(dx)^{-\frac{2}{3}})\psi(dx)^{-\frac{2}{3}}+&\\[2mm]
\epsilon_3\,\mathrm{Gz}(\psi\,(dx)^{-\frac{2}{3}},\chi\,(dx)^{-\frac{2}{3}})\varphi\,
(dx)^{-\frac{2}{3}}& \text{for }(\gamma,\lambda,\tau)
=(-\frac{2}{3},-\frac{2}{3},-\frac{2}{3}), \\[3mm]
\epsilon_1  \{ d\{\varphi(dx)^{-\frac{1}{2}},\psi(dx)^{-\frac{1}{2}}
\},\chi (dx)^{-\frac{1}{2}}\}+&\\[2mm]
\epsilon_2 \{d
\{\chi(dx)^{-\frac{1}{2}},\varphi (dx)^{-\frac{1}{2}}\},\psi(dx)^{-\frac{1}{2}} \}+&\\[2mm]
\epsilon_3\{d \{\psi(dx)^{-\frac{1}{2}},\chi(dx)^{-\frac{1}{2}}
\},\varphi(dx)^{-\frac{1}{2}}  \}& \text{for
}(\gamma,\lambda,\tau)=(-\frac{1}{2},-\frac{1}{2},-\frac{1}{2}),\\[3mm]
\epsilon_1 \{ d\{\varphi(dx)^\lambda,\chi(dx)^\tau\},\psi(dx)^\gamma
\}+&\\[2mm]
\epsilon_2 \{ d \{
\psi(dx)^\gamma,\chi(dx)^\tau\},\varphi(dx)^\lambda \}& \text{for
}1+\lambda+\tau=0\text{ and } \gamma=\lambda\not =\tau, \\[3mm]
\Delta_{\lambda, \gamma, \tau;3}(
\varphi(dx)^\lambda,\psi(dx)^\gamma,\chi(dx)^\tau)
&\text{otherwise.}
\end{cases}
\end{array}
$$
\end{thm}

{\bf Proof.} The $\mathfrak{vect}(1)$-invariance of the differential
operator
\[A=\sum_{i+j+l=3}\alpha_{i,j,l}\,\varphi^{(i)}\psi^{(j)}\chi^{(l)}\,
(dx)^{\mu}\] is equivalent to invariance with respect to just the
two fields that generate $\mathfrak{vect}(1)$, namely $\frac{d}{dx}$
and $x^3\frac{d}{dx}$, and is equivalent to the system:

\begin{equation}\label{sys}\tiny
\begin{array}{l}
 \alpha_{i,j,l} \left (\lambda \binom{i}{r-1}+
\binom{i}{r}\right) + \alpha_{i-r+1,j+r-1,l} \left (\gamma
\binom{j+r-1}{r-1}+ \binom{j+r-1}{r}\right) +\alpha_{i-r+1,j,l+r-1}
\left (\tau \binom{l+r-1}{r-1}+ \binom{l+r-1}{r}\right)=0,\\[3mm]
\normalfont\text{where $i+j+l=3$ and $r=2,3$.}
\end{array}
\end{equation}

Our strategy: We find all solutions to the system above with all
possible values of the parameters $\lambda, \gamma$ and $\tau$;
assume that the space of solutions is $m$-dimensional. Then we try
to construct $m$ linearly independent operators expressed in terms
of the de Rham differential, the Poisson bracket, the Grozman
operator or the Feigin-Fuchs anti-symmetric operators: this is only
possible for particular values of $\lambda, \mu$ and $\tau$. If this
is impossible (we test ALL possible combination of the exterior
differential and the bilinear operators), or the space obtained is
of dimension $<m$, the operator $A$ we found is new. These values
are as follows:

1. The case where $\lambda=\gamma=\tau=0$. Here the system above
admits the following solutions
\[
\alpha_{0,0,3}=\alpha _{3,0,0}=\alpha_{0,3,0}=0, \quad \alpha
_{0,2,1}=-\alpha _{0,1,2}, \quad \alpha _{2,0,1}=-\alpha _{1,0,2},
\quad \alpha_{2,1,0}=-\alpha _{1,2,0}.
\]
We expect then a family of operators that depends on
$\alpha_{1,2,0}$, $\alpha_{1,0,2}$, $\alpha_{1,1,1}$ and $
\alpha_{0,1,2}$. On the other hand, the following operators are
invariant and cannot be multiples of each other:
\begin{equation}\label{invcomp}
\left \{d\varphi,d\psi\right \}\chi, \quad \left
\{d\chi,d\varphi\right \}\psi , \quad \left \{d\psi,d\chi \right
\}\varphi, \quad d\varphi\otimes d\psi\otimes d\chi.
\end{equation}
Therefore, $A$ is a linear combination of the invariant operators
(\ref{invcomp}).

2. The case where $\lambda=\gamma=0$ and $\tau\not=0$. (The cases
$\lambda=\tau=0, \gamma\not=0$ and $\gamma=\tau=0, \lambda\not=0$
follow from permutations. The case $\tau=-2$ should be disregarded
since it is the dual to case 1.) Without going into details, observe
that the system (\ref{sys}) has solutions depending on two
parameters. This means that $A$ is a linear combination of the
invariant operators
\[
 \{d\varphi, \chi(dx)^\tau\} d \psi \text{ and }
\{d\psi,\chi(dx)^\tau\} d\varphi.
\]

3. The case where $\lambda=0$ and $\gamma, \tau\not=0$. We proceed
as in case 2. 

4. The case where $\lambda,\gamma, \tau\not=0$. We proceed as in
case 2. \qed 

\begin{rmk}
{\rm The reader may wonder why the operators $\Delta_{\l,3}$ for
$\lambda=-\frac{1}{2}$, $-\frac{2}{3}$, and $0$ do not appear in our
list. But in our classification some operators generalize the
operators $\Delta_{\l,3}$ for particular values of $\epsilon$. For
instance, $\Delta_{-\frac{2}{3},3}(\varphi, \psi,\chi)$ coincides
with the operator $$\epsilon_1\,\mathrm{
Gz}(\varphi,\psi)\chi+\epsilon_2\,\mathrm{Gz}(\chi,\varphi)\psi
+\epsilon_3\,\mathrm{Gz}(\psi,\chi)\varphi\text{~~ for
$\epsilon_1=\epsilon_2=\epsilon_3$.}$$ }
\end{rmk}

\subsubsection{Order 4}
The following ternary invariant differential operator is new, and so
are its dualizations and permutations (here $\mu=\tau+\frac{5}{2}$
and $\tau\not= -\frac{3}{4})$:
\[\small
\Xi: \cF_{-\tau-\frac{3}{2}}\otimes \cF_{\tau}\otimes
\cF_{\tau}\rightarrow \cF_{\mu} \quad
\varphi(dx)^{-\tau-\frac{3}{2}}\otimes \psi(dx)^\tau \otimes \chi
(dx)^\tau \mapsto
\!\!\sum_{i+j+l=4}\!\!\!\!\alpha_{i,j,l}\,\varphi^{(i)}\psi^{(j)}\chi^{(l)}
(dx)^{\mu},\] where the coefficients are given by
\[
\begin{array}{lcllcllcl}
\alpha_{4,0,0}&=&0,&
\alpha_{0,4,0}&=&-\tau(3+2\tau),&\alpha_{1,1,2}&=&-10(1+\tau),\\
\alpha_{0,0,4}&=&\tau (3+2\tau), & \alpha_{3,1,0}&=&
-\frac{8}{3}\tau(2+3\tau),
&\alpha_{3,0,1}&=&\frac{8}{3}\tau(2+3\tau),\\[2mm]
\alpha_{1,3,0}&=&-\frac{2}{3}\tau(13+12\tau),&
\alpha_{1,0,3}&=&\frac{2}{3}\tau(13+12\tau),&
\alpha_{0,3,1}&=&\frac{5}{3}(3+2\tau),\\[2mm]
\alpha_{0,1,3}&=&-\frac{5}{3}(3+2\tau),&
\alpha_{2,2,0}&=&-4\tau(2+3\tau),
&\alpha_{2,0,2}&=&4\tau(2+3\tau), \\[2mm]

\alpha_{0,2,2}&=&0,& \alpha_{2,1,1}&=&0,&
\alpha_{1,2,1}&=&10(1+\tau).\\[2mm]
\end{array}
\]
For $\tau=-\frac{3}{4}$, we have a 2-parametric family of invariant
operators:
\[\small
\begin{array}{lcllcllcl}
\alpha_{0,0,4}&=& - s- t,&\alpha_{3,1,0}&=&\frac{4}{9}
   \left(4 t- s\right),& \alpha_{3,0,1}&=&\frac{4}{9} \left(
   s+5 t\right),\\[3mm]
\alpha_{1,3,0}&=&\frac{4}{9} \left(4 s- t\right),&
\alpha_{1,0,3}&=&-\frac{4}{9} \left(4
s+5t\right),&\alpha_{0,3,1}&=&\frac{4}{9}
\left(5 s+ t\right),\\[3mm]
\alpha_{0,1,3}&=&-\frac{4}{9}\left(5 s+4 t\right),&
\alpha_{2,2,0}&=&-\frac{2}{3} \left(s+t\right),&
\alpha_{2,0,2}&=&\frac{2}{3}
   s,\\[3mm]
\alpha_{0,2,2}&=&\frac{2}{3} t, &  \alpha _{2,1,1}&=&\frac{20}{9}
t,&
\alpha_{1,2,1}&=&\frac{20}{9} s,\\[3mm]
\alpha_{1,1,2}&=&-\frac{10}{9} \left(s+t\right),&\alpha_{0,4,0}&=&s,
&\alpha_{4,0,0}&=&t.
   \end{array}
\]

\begin{thm} Every order $4$ invariant differential operator
is, up to permutations and dualizations, a multiple of the following
operators: $\varphi (dx)^\lambda\otimes \psi (dx)^\gamma\otimes \chi
(dx)^\tau\mapsto$
$$
\left\{
\begin{array}{ll}
\epsilon_1 \{\{d\varphi,d\psi \},\chi(dx)^\tau \}+\epsilon_2 \{
 \{\chi(dx)^\tau,d\varphi\},d\psi  \}
& \text{ for }(\gamma,\lambda)=(0,0), \\[3mm]

\epsilon_1 \{d \{d\psi,\varphi(dx)^{-2}\},\chi
(dx)^{-2}\}+&\\[2mm]

\epsilon_2 \{ d \{d\psi,\chi (dx)^{-2} \},\varphi
(dx)^{-2}  \} & \text{ for }(\lambda,\gamma,\tau)=(-2,0,-2), \\[3mm]

\Delta_{\lambda, \gamma, \tau;3}(\varphi(dx)^\lambda, d\psi,
\chi(dx)^\tau) & \text{ for }\gamma=0,\lambda,
\tau\not=0,\\[3mm]

 \{{\mathrm Gz} (
\varphi(dx)^{-\frac{2}{3}},\psi(dx)^{-\frac{2}{3}}),\chi(dx)^\tau \}
& \text{ for
}(\lambda,\gamma) =(-\frac{2}{3},-\frac{2}{3}) \text{ and } \tau \not =-1,\\[3mm]

\epsilon_1 {\mathrm Gz} ( \{
\varphi(dx)^{-1},\psi(dx)^{-\frac{2}{3}} \},\chi(dx)^{-\frac{2}{3}})+&\\[2mm]
\epsilon_2 {\mathrm Gz} ( \{
\varphi(dx)^{-1},\chi(dx)^{-\frac{2}{3}}
\},\psi(dx)^{-\frac{2}{3}})& \text{ for }(\lambda,\gamma,\tau)
=(-1,-\frac{2}{3},-\frac{2}{3}), \\[3mm]

\Xi (\varphi(dx)^\lambda,\psi(dx)^\gamma,\chi(dx)^\tau) & \text{ for
}(\lambda,\gamma)=(-\frac{3}{2}-\tau,\tau)\\[2mm]
& \text{ and }\tau\not=-\frac{3}{2},-\frac{2}{3},0,-\frac{3}{4},\\[3mm]

\Xi_{s,t}(\varphi(dx)^\lambda,\psi(dx)^\gamma,\chi(dx)^\tau) &
\text{ for }(\lambda,\gamma,\tau)=(-\frac{3}{4} , -\frac{3}{4} ,
-\frac{3}{4} ).
\end{array}
\right.
$$
\end{thm}
\subsubsection{Order 5}
The following order 5 ternary invariant differential operator is
new, and so are its dualizations and permutations (unless they
coincide):
\begin{equation}
\label{fifth}\small \Gamma:\cF_{-\frac{2}{3}}\otimes
\cF_{-\frac{2}{3}}\otimes \cF_{-\frac{4}{3}}\rightarrow
\cF_{\frac{7}{3}} \quad \varphi(dx)^{-\frac{2}{3}}\otimes
\psi(dx)^{-\frac{2}{3}} \otimes \chi(dx)^{-\frac{4}{3}} \mapsto
\sum_{i+j+l=5}\!\alpha_{i,j,l}\varphi^{(i)}\psi^{(j)}\chi^{(l)}
(dx)^{\frac{7}{3}},
\end{equation}
where the constants are given by
\[
\begin{array}{lcrccrccrccrccr}
\alpha_{0,0,5}&=&\frac{2}{5},&\alpha_{0,2,3}&=&-1,&
\alpha_{0,2,3}&=&-\frac{5}{2},&
\alpha_{0,4,1}&=&-\frac{17}{10},&\alpha_{0,5,0}&=&-\frac{2}{5},\\[3mm]
\alpha_{1,0,4}&=&1,&\alpha_{1,1,3}&=&\frac{3}{2},
&\alpha_{1,2,2}&=&-\frac{9}{4},& \alpha_{1,3,1}&=&-\frac{9}{4},&
\alpha_{1,4,0}&=&-\frac{3}{5},\\[3mm]

\alpha_{2,0,3}&=&-1,&\alpha_{2,1,2}&=&-\frac{9}{4},&
\alpha_{2,2,1}&=&\frac{9}{2},& \alpha_{2,3,0}&=&3,&
\alpha_{3,0,2}&=&-\frac{5}{2},\\[3mm]

\alpha_{3,1,1}&=&-\frac{9}{4},& \alpha_{3,2,0}&=&3,&
\alpha_{4,0,1}&=&-\frac{17}{10},&\alpha_{4,1,0}&=&-\frac{3}{5},&
\alpha_{5,0,0}&=&-\frac{2}{5}.
\end{array}
\]
\begin{thm}
Every order $5$ invariant differential operator is, up to
permutations and dualizations, a multiple of the following
operators: $\varphi(dx)^\lambda\otimes \psi(dx)^\gamma\otimes \chi
(dx)^\tau\mapsto $
$$
\left \{
\begin{array}{ll}
\epsilon_1\left \{\left \{d\varphi,d\psi\right \},d\chi\right
\}+\epsilon_2\left \{ \left \{d\chi,d\varphi\right \},d\psi \right
\}
& \text{for }(\lambda,\gamma,\tau)=(0,0,0), \\[3mm]
\epsilon_1 \left \{d\left \{d\varphi,\chi(dx)^{-2}\right
\},d\psi\right \}+\\[2mm]\epsilon_2 \left \{ d\left
\{d\psi,\chi(dx)^{-2}\right \},d\varphi \right \}
& \text{for }(\lambda,\gamma,\tau)=(0,0,-2), \\[2mm]
\Delta_{\lambda, \gamma, \tau;3}(d\varphi, d\psi, \chi(dx)^{\tau}) &
\text{for }(\lambda,\gamma)=(0,0)
 \text{ and }\tau\not= -4,-2,0, \\[3mm]
 \{{\mathrm Gz} (
\varphi(dx)^{-\frac{2}{3}},\psi(dx)^{-\frac{2}{3}}),d\chi \} &
\text{for
}(\lambda,\gamma,\tau)=(-\frac{2}{3},-\frac{2}{3},0),\\[3mm]
\Xi (\varphi(dx)^{-\frac{5}{2}},d\psi,\chi dx) & \text{for
}(\lambda,\gamma,\tau)=(-\frac{5}{2},0,1),\\[3mm]
\Gamma (\varphi(dx)^\lambda,\psi (dx)^\gamma,\chi(dx)^\tau) &
\text{for
}(\lambda,\gamma,\tau)=(-\frac{2}{3},-\frac{2}{3},-\frac{4}{3}),\\[3mm]
\Theta_{\pm} (\varphi(dx)^\lambda,\psi (dx)^\gamma,\chi(dx)^\tau) &
\text{for }
\lambda=\gamma=\tau=-\frac{1}{12}(9\pm\sqrt{21}).\\[3mm]
\end{array}
\right.
$$
\end{thm}
\subsubsection{Order 6}
\begin{thm} Every order $6$ invariant differential operator is, up to permutations
and dualizations, a multiple of the following operators:
$\varphi(dx)^\lambda\otimes \psi(dx)^\gamma\otimes
\chi(dx)^\tau\mapsto$
$$
\left \{
\begin{array}{ll}
\Delta_{1,3}(d \varphi, d \psi, d \chi) & \text{ for
}(\lambda,\gamma,\tau)=(0,0,0), \\[2mm]
\Xi (d \varphi,d \psi,\chi(dx)^{-\frac{5}{2}})& \text{ for
}(\lambda,\gamma,\tau)=(0,0,-\frac{5}{2}),\\[2mm]
\Upsilon
(\varphi(dx)^{-\frac{5}{4}},\psi(dx)^{-\frac{5}{4}},\chi(dx)^{-\frac{5}{4}})
& \text{ for
}(\lambda,\gamma,\tau)=(-\frac{5}{4},-\frac{5}{4},-\frac{5}{4}).\\[2mm]
\end{array}
\right .
$$
\end{thm}
\subsubsection{Order $>6$}
There are no such ternary $\mathfrak{vect}(1)$-invariant
differential operators.

{\bf Proof.} This statement can be checked by a direct computation
for operators of order $k=7$. To prove the result for any order
$k>7$, we proceed as follows.

{\bf Step 1}. The case where $\lambda\gamma\tau\not =0$. The
invariance with respect to the vector field $x^{k+1}\frac{d}{dx}$,
where $k$ is the order of the operator, is equivalent to the system:
\begin{equation}\label{sys2}\tiny
\begin{array}{l}
 \alpha_{i,j,l} \left (\lambda \binom{i}{r-1}+
\binom{i}{r}\right) + \alpha_{i-r+1,j+r-1,l} \left (\gamma
\binom{j+r-1}{r-1}+ \binom{j+r-1}{r}\right) +\alpha_{i-r+1,j,l+r-1}
\left (\tau \binom{l+r-1}{r-1}+ \binom{l+r-1}{r}\right)=0,\\[3mm]
\normalfont\text{where $i+j+l=k$ and $r=2,\ldots, k+1$.}
\end{array}
\end{equation}
Let us first prove that
\begin{equation}\label{k7}
a_{k,0,0}=a_{0,k,0}=a_{0,0,k}=a_{k-1,i,j}=a_{i,k-1,j}=
a_{i,j,k-1}=a_{k-2,i,j}=a_{i,k-2,j}=a_{i,j,k-2}=0.
\end{equation}
To do so, consider a sub-system of (\ref{sys2}) whose coefficients
are as in (\ref{k7}) for $r=k+1,k,k-1$ and $r=3,2$. We get then a
linear system with a $23\times 18$-matrix (too huge to be inserted
here).  When $\lambda\gamma\tau\not =0$, this matrix is of rank 18,
and hence all solutions are zero as in (\ref{k7}).

Let
\[
\alpha_{k-t, i, t-i}=\alpha_{i, k-t, t-i}=\alpha_{i, t-i, k-t}=0
\text{ for all $t,i=0,\ldots,n$}.
\]
Let us induct on $n$ to prove that the remaining coefficients are
zero.

\underline{For $r=k-n$: } Take $i=k-n-1$ and any $j=s \not=0,n+1$ in
the system (\ref{sys2}). We get:
\[
\lambda\,  \alpha_{k-n-1, s, n+1-s} + \left ( \gamma
\binom{s+k-n-1}{k-n-1}+ \binom{s+k-n-1}{k-n} \right ) \alpha_{0,
k-n-1+s, n+1-s}=0.
\]
The induction hypothesis implies that $\alpha_{0, k-n-1+s,
n+1-s}=0$, and hence
\begin{equation}\label{zer1}
 \alpha_{k-n-1, s, n+1-s}=0 \text{ for all } s
\not=0,n+1.
\end{equation}
Similarly, we can prove that $\alpha_{s, k-n-1,
n+1-s}=\alpha_{s,n+1-s,k-n-1}=0 \text{ for all } s \not=0,n+1.$

\underline{For $r=3$: } Take $i=2$ and $j=n-1$ in the system
(\ref{sys2}). We get:
\begin{equation}\label{sys3}
\begin{array}{clcclc}
\left (\gamma \binom{n+1}{2}+\binom{n+1}{3}\right )\, \alpha_{0,
n+1, k-n-1}&=&0, &\left (\gamma \binom{n+1}{2}+\binom{n+1}{3}\right
)\, \alpha_{k-n-1, n+1, 0}&=&0,\\[2mm]

\left (\lambda \binom{n+1}{2}+\binom{n+1}{3}\right )\, \alpha_{
n+1,0, k-n-1}&=&0, &\left (\lambda
\binom{n+1}{2}+\binom{n+1}{3}\right
)\, \alpha_{n+1, k-n-1, 0}&=&0,\\[2mm]

\left (\tau \binom{n+1}{2}+\binom{n+1}{3}\right )\,
\alpha_{0,k-n-1,n+1}&=&0, &\left (\tau
\binom{n+1}{2}+\binom{n+1}{3}\right )\, \alpha_{k-n-1, 0,n+1}&=&0.
\end{array}
\end{equation}

\underline{For $r=2$: } Take $i=1$ and $j=n$ in the system
(\ref{sys2}). We get:
\begin{equation}\label{sys4}
\begin{array}{clcclc}
\left (\gamma (n+1)+\binom{n+1}{2}\right )\, \alpha_{0, n+1,
k-n-1}&=&0, &\left (\gamma (n+1)+\binom{n+1}{2}\right
)\, \alpha_{k-n-1, n+1, 0}&=&0,\\[2mm]

\left (\lambda (n+1)+\binom{n+1}{2}\right )\, \alpha_{ n+1,0,
k-n-1}&=&0, &\left (\lambda (n+1)+\binom{n+1}{2}\right
)\, \alpha_{n+1, k-n-1, 0}&=&0,\\[2mm]

\left (\tau (n+1)+\binom{n+1}{2}\right )\, \alpha_{0,k-n-1,n+1}&=&0,
&\left (\tau (n+1)+\binom{n+1}{2}\right )\, \alpha_{k-n-1,
0,n+1}&=&0.
\end{array}
\end{equation}
Eqs. (\ref{sys3}) and Eqs. (\ref{sys4}) imply that
\begin{equation}\label{zer2}\tiny
\alpha_{k-n-1,0,n+1}=\alpha_{k-n-1,n+1,0}=\alpha_{0,k-n-1,n+1}=
\alpha_{n+1,k-n-1,0}= \alpha_{0,n+1,k-n-1}=\alpha_{n+1,0,k-n-1}=0.
\end{equation}
Form Eqs. (\ref{zer1}) and Eqs. (\ref{zer2}) we see that the
induction hypothesis is true at $n+1$, hence the result.

{\bf Step 2}. The case where $\lambda\gamma\tau\not =0$. We provide
a proof here only when $\tau=0$ and $\lambda\gamma\not =0$. Now the
operator $A(\cdot,\cdot,1)$ turns into a
$\mathfrak{vect}(1)$-invariant binary one since $\tau=0$. It follows
from Grozman's classification that $A(\cdot,\cdot,1)$  is
identically zero because $k>7$. Therefore, $\alpha_{i,k-i,0}=0$ for
$i=0,\ldots, k$. Thus, $A=B(\cdot,\cdot,d)$, where $B$ is an order
$k-1$ operator defined on ${\mathcal F}_\lambda\otimes {\mathcal
F}_\gamma \otimes {\mathcal F}_1$. The
$\mathfrak{vect}(1)$-invariance of $A$ automatically implies
$\mathfrak{vect}(1)$-invariance of $B$. Thanks to Step 1, the
operator $B$ must be identically zero, together with $A$. \qed

\section{Ternary invariant differential operators on
multi-dimensional manifolds}
Consider now $M=\mathbb{R}^n$. Denote by $g$ the usual
pseudo-Riemannian metric on $\mathbb{R}^n$ of signature $(p,q)$,
where $p+q=n$. The conformal transformations are generated by vector
fields \footnote{We use the conventional summation over repeated
indices. Indices are raised or lowered by means of the metric $g$.}:
\begin{equation}
\label{inf}
\begin{array}{lcllcl}
X_i&=&\displaystyle \frac{\partial}{\partial x^i},& X_{ij}&=&
\displaystyle x_i\frac{\partial}{\partial x^j}-
x_j\frac{\partial}{\partial x^i},\\[2mm]
X_0&=& \displaystyle x^i\frac{\partial}{\partial x^i},&
\overline{X_i}&=&\displaystyle x_jx^j\frac{\partial}{\partial
x^i}-2x_ix^j \frac{\partial}{\partial x^j},
\end{array}
\end{equation}
where $(x^1,\ldots,x^n)$ are coordinates on $\mathbb{R}^n$ and
$x_i=g_{is}x^s$. In what follows the Lie algebra
$\mathfrak{o}(p+1,q+1)$ is supposed to be realized by (\ref{inf}).
Obviously, the Lie algebra generated by the vector fields $X_{ij}$
is isomorphic to $\mathfrak{o}(p,q)$. We have therefore an inclusion
of Lie algebras:
\[
\mathfrak{o}(p,q)\subset \mathfrak{o}(p+1,q+1)\subset
\mathfrak{vect}(n).
\]
We first classify all ternary $\mathfrak{o}(p+1,q+1)$-invariant
differential operators, then we use the result to classify those
that are $\mathfrak{vect}(1)$-invariant.
\subsection{Conformally invariant operators}
The following technique is due to Ovsienko and Redoux \cite{or}. We
identify ternary differential operators with their symbols. This
identification allows us to study the subalgebra $I$ of ${\mathfrak
o}(p,q)$-invariant polynomials
\[
I:=\mathbb{C}[x^1,\ldots,x^n,\xi_1,\ldots,\xi_n,\eta_1,\ldots,\eta_n,
\zeta_1,\ldots,\zeta_n]^{{\mathfrak o}(p,q)}
\]
and to apply the Weyl invariant theory \cite{wey}. It follows that
$I$ has the following generators:
\begin{equation}\label{gen}
\begin{array}{lcllcllcllcllcl}
{\mathrm R}_{xx}&=&x^ix_j, &{\mathrm R}_{x\xi}&=&x^i\xi_i,& {\mathrm
R}_{x\eta}&=&x^i\eta_i,& {\mathrm R}_{x\zeta}&=&x^i\zeta_i,&
{\mathrm R}_{\xi\xi}&=&\xi^i\xi_j,\\[2mm]
{\mathrm R}_{\xi\eta}&=&\xi^i\eta_i,& {\mathrm
R}_{\xi\zeta}&=&\xi^i\zeta_i,& {\mathrm
R}_{\eta\eta}&=&\eta^i\eta_i,& {\mathrm
R}_{\eta\zeta}&=&\eta^i\zeta_i,& {\mathrm
R}_{\zeta\zeta}&=&\zeta^i\zeta_i.
\end{array}
\end{equation} The $X_i$-invariance implies that
$\frac{\partial P}{\partial x^i}=0$ for every polynomial $P$.
Therefore, the generators from the list (\ref{gen}) that contain $x$
should be disregarded. Hence, any $\mathfrak{o}(p+1,q+1)$-invariant
differential operator should be of form
\begin{equation}
\label{conf} B=\sum_{a,b,c,d,e,f\geq 0} \alpha_{abcdef} {\mathrm
R}^{a,b,c,d,e,f}
\end{equation}
where ${\mathrm R}^{a,b,c,d,e,f}={\mathrm R}_{\xi\xi}^a\, {\mathrm
R}_{\xi\eta}^b\, {\mathrm R}_{\xi\zeta}^c\,{\mathrm
R}_{\eta\eta}^d\, {\mathrm R}_{\eta\zeta}^e\, {\mathrm
R}_{\zeta\zeta}^f$.
\begin{thm}
\label{confo} For almost every $\lambda,\gamma$ and $\tau$, there
exist ${\mathfrak o}(p+1,q+1)$-invariant operators:
\begin{equation}\label{B2k}
B_{2k}:\cF_\l\otimes\cF_\gamma\otimes \cF_\tau\rightarrow
\cF_{\l+\gamma+\tau+\frac{2k}{n}}  \quad B_{2k}=\sum_{a+b+c+d+e+f=k}
\alpha_{abcdef} {\mathrm R}^{a,b,c,d,e,f}
\end{equation}
where $k=0,1,\ldots$ and the constants are given by the induction
formulas:
\begin{equation}\label{NEqs}
\small{\begin{array}{rcl} 2(a+1)(2(a+1)+n(2\lambda-1))\,
\alpha_{a+1,b,c,d,e,f}-(b+2)(b+1)\,
\alpha_{a,b+2,c,d-1,e,f}&& \\[2mm]
-2(b+1)(c+1) \,\alpha_{a,b+1,c+1,d,e-1,f}+2(b+1)(b+e+2d+n\gamma)\,
\alpha_{a,b+1,c,d,e,f}&&\\[2mm]
-(c+2)(c+1) \,\alpha_{a,b,c+2,d,e,f-1}+ 2(c+1)(c+e+2f+n\tau)\,
\alpha_{a,b,c+1,d,e,f}&=&0,\\[2mm]
2(d+1)(2(d+1)+n(2\gamma-1))\, \alpha_{a,b,c,d+1,e,f}-(b+2)(b+1)\,
\alpha_{a-1,b+2,c,d,e,f}&& \\[2mm]
-2(b+1)(e+1) \,\alpha_{a,b+1,c-1,d,e+1,f}+2(b+1)(b+c+2a+n\lambda)\,
\alpha_{a,b+1,c,d,e,f}&&\\[2mm]
-(e+2)(e+1) \,\alpha_{a,b,c,d,e+2,f-1}+ 2(e+1)(c+e+2f+n\tau)\,
\alpha_{a,b,c,d,e+1,f}&=&0,\\[2mm]
2(f+1)(2(f+1)+n(2\tau-1))\, \alpha_{a,b,c,d,e,f+1}-(e+2)(e+1)\,
\alpha_{a,b,c,d-1,e+2,f}&& \\[2mm]
-2(e+1)(c+1) \,\alpha_{a,b-1,c+1,d,e+1,f}+2(e+1)(b+e+2d+n\gamma)\,
\alpha_{a,b,c,d,e+1,f}&&\\[2mm]
-(c+2)(c+1) \,\alpha_{a-1,b,c+2,d,e,f}+ 2(c+1)(c+b+2a+n\lambda)\,
\alpha_{a,b,c+1,d,e,f}&=&0.
\end{array}}
\end{equation}
\end{thm}
{\bf Example.} For $k=1$, we have (here $s,t$ and $u$ are
parameters):
\[
\begin{array}{ccl}
B_2&=&n(\gamma s+\tau t)(2+n(2\gamma-1))(2+n(2\tau-1))
{\mathrm R}_{\xi\xi}\\[2mm]
&&+n(\lambda s+\tau u)(2+n(2\lambda-1))(2+n(2\tau-1)){\mathrm
R}_{\eta\eta}\\[2mm]
&&+n(\gamma u+\lambda t)(2+n(2\lambda-1))(2+n(2\gamma-1)){\mathrm
R}_{\zeta\zeta}\\[2mm]
&&-(2+n(2\lambda-1))(2+n(2\gamma-1))(2+n(2\tau-1))(s {\mathrm
R}_{\xi\eta}+t{\mathrm R}_{\xi\zeta}+u{\mathrm R}_{\eta\zeta}).
\end{array}
\]
{\bf Proof of Theorem \ref{confo}.} The $X_0$-invariance implies
that ($L_X$ is the Lie derivative along the field $X$)
\[
L_{X_0}B=\sum_{a,b,c,d,e,f\geq 0}
(n(\mu-\l-\gamma-\tau)-2(a+b+c+d+e+f))\,c_{a,b,c,d,e,f} {\mathrm
R}^{a,b,c,d,e,f}=0.
\]
Hence $(n(\mu-\l-\gamma-\tau)-2(a+b+c+d+e+f))=0$. Thus, the operator
is homogeneous (all its components are of the same degree). In order
to deduce the induction formulas, we need the following Proposition
(see, e.g., \cite{dlo}):
\begin{pro}\label{kan}
The action of $\overline{X_i}$ reads as follows:
\[
L_{\overline{X_i}}^{\l,\gamma,\tau;\mu}=
L_{\overline{X_i}}^{\mu-\l-\gamma-\tau}-\xi_i T_\xi -\eta_i T_\eta
-\zeta_iT_\zeta+2((E_\xi+n\l)\partial_{\xi^i}+(E_\eta+n\gamma)\partial_{\eta^i}+(
E_\zeta+n\tau)\partial_{\zeta^i}),
\]
where
\[
\begin{array}{ccl}
L_{\overline{X_i}}^{\mu-\l-\gamma-\tau}&=&x_j
x^j\partial_i-2x_ix^j\partial_j-2n(\mu-\l-\gamma-\tau) x_i\\[2mm]
&& -2 \left ((\xi_i x_j -\xi_j x_i)\partial_{\xi_j} +(\eta_i x_j
-\eta_j x_i)\partial_{\eta_j}+(\zeta_i x_j -\zeta_j
x_i)\partial_{\zeta_j}
\right )\\[2mm]
&& +2\left (\xi_j x^j \partial_{\xi^i}+\eta_j x^j
\partial_{\eta^i}+\zeta_j x^j \partial_{\zeta^i}\right ),
\end{array}
\]
and where $T_{\xi}=\partial_{\xi^j}\partial_{\xi_j}$ is the trace
and $E_{\xi}=\xi_j \partial_{\xi_j}$ is the Euler operator.
\end{pro}
Proposition \ref{kan} implies that 
$L_X^{\lambda,\gamma,\tau;\mu}{\mathrm R}^{a,b,c,d,e,f}=$
\begin{equation}\label{L_x}
\footnotesize\begin{array}{l} 2(2k-n(\mu-\l-\gamma-\tau)) {\mathrm
R}^{a,b,c,d,e,f}\,x_i+\\[2mm]
\left ( 2a(2a+n(2\lambda-1))\, {\mathrm R}^{a-1,b,c,d,e,f}-b(b-1)\,
{\mathrm
R}^{a,b-2,c,d+1,e,f}\right.\\[2mm]
-2bc \, {\mathrm R}^{a,b-1,c-1,d,e+1,f}+2b(b-1+e+2d+n\gamma)\,
{\mathrm R}^{a,b-1,c,d,e,f}\\[2mm]
\left.  -c(c-1) \, {\mathrm R}^{a,b,c-2,d,e,f+1}+
2c(c-1+e+2f+n\tau)\, {\mathrm R}^{a,b,c-1,d,e,f}
\right )\xi_i+\\[2mm]
\left ( 2d(2d+n(2\gamma-1))\, {\mathrm R}^{a,b,c,d-1,e,f}-b(b-1)\,
{\mathrm
R}^{a+1,b-2,c,d,e,f}\right.\\[2mm]
-2be \, {\mathrm R}^{a,b-1,c+1,d,e-1,f}+2b(b-1+c+2a+n\lambda)\,
{\mathrm R}^{a,b-1,c,d,e,f}\\[2mm]
\left.  -e(e-1) \, {\mathrm
R}^{a,b,c,d,e-2,f+1}+2e(e-1+c+2f+n\tau)\,
{\mathrm R}^{a,b,c,d,e-1,f}\right )\eta_i+\\[2mm]
\left ( 2f(2f+n(2\tau-1))\, {\mathrm R}^{a,b,c,d,e,f-1}-e(e-1)\,
{\mathrm
R}^{a,b,c,d+1,e-2,f}\right.\\[2mm]
-2ec \, {\mathrm R}^{a,b+1,c-1,d,e-1,f}+2e(e-1+b+2d+n\gamma)\,
{\mathrm R}^{a,b,c,d,e-1,f}\\[2mm]
\left.  -c(c-1) \, {\mathrm
R}^{a+1,b,c-2,d,e,f}+2c(c-1+b+2a+n\lambda)\,
{\mathrm R}^{a,b,c-1,d,e,f}\right )\zeta_i.\\[2mm]
\end{array}
\end{equation} Now acting by $L_X^{\lambda,\gamma,\tau;\mu}$ on
$B_{2k}$ we get, upon collecting terms, a recurrent system of linear
equations (\ref{NEqs}). 
The operators $B_{2k}$ depend on some parameters and may vanish for
some particular values of $\lambda, \gamma$ and $\tau$.
\subsection{$\mathfrak{vect}(n)$-invariant differential operators}
As every $\mathfrak{vect}(n)$-invariant operator is
$\mathfrak{o}(p+1,q+1)$-invariant, it suffices to check whether the
conformally invariant operators listed in Theorem \ref{confo} are
$\mathfrak{vect}(n)$-invariant or not:
\begin{thm}
The only $\mathfrak{vect}(n)$-invariant ternary differential
operator is the multiplication operator:
\[
\cF_\l\otimes \cF_\gamma\otimes \cF_\tau \rightarrow
\cF_{\lambda+\gamma+\tau} \quad \varphi
|\mathrm{vol}|^\lambda\otimes \psi |\mathrm{vol}|^\gamma \otimes
\chi |\mathrm{vol}|^\tau \mapsto \varphi\cdot\psi\cdot\chi
|\mathrm{vol}|^{\lambda+\gamma+\tau}.
\]
\end{thm}
{\bf Proof.} 
A direct computation shows that the action by
$X\in\mathfrak{vect}(n)$ is as follows (cf. with the action in
Proposition \ref{kan}):
\[
\begin{array}{lcl}
L_{X}^{\lambda,\gamma,\tau;\mu}&=&(\mu-\lambda-\gamma-\tau)\,
\mathrm{Div}(X) -\left (\partial_j X^m \xi_m \partial_{\xi_j}
+\partial_j X^m \eta_m
\partial_{\eta_j}+\partial_j X^m \zeta_m
\partial_{\zeta_j} \right )\\[2mm]
&&+(\text{higher order derivatives } \partial_{i_1}\cdots
\partial_{i_l} X).
\end{array}
\]
Now, since every $\mathfrak{vect}(n)$-invariant operator is
$\mathfrak{o}(p+1,q+1)$-invariant, we have (see (\ref{B2k})):
\[
L_X^{\lambda, \gamma, \tau}(B_{2k})=\sum_{a+b+c+d+e+f=k}
(\mu-\lambda-\gamma-\tau)
\alpha_{a,b,c,d,e,f}\mathrm{Div}(X)\mathrm{R}^{a,b,c,d,e,f}+
(\text{other terms}).
\]
Thus $(\mu-\lambda-\gamma-\tau) \alpha_{a,b,c,d,e,f}=0$. Since
$\mu-\lambda-\gamma-\tau=\frac{2k}{n}$, it follows that either $k=0$
or $\alpha_{a,b,c,d,e,f}=0$. Thus if $B_{2k}$ is not the scalar
operator (\ref{scal}), it must be identically zero.\qed


\end{document}